\documentclass[12pt, reqno]{amsart}
\usepackage{amsmath, amsthm, amscd, amsfonts, amssymb, graphicx, color}
\usepackage[bookmarksnumbered, colorlinks, plainpages]{hyperref}
\input{mathrsfs.sty}

\theoremstyle{definition}

\theoremstyle{remark}

\numberwithin{equation}{section}

\begin{document}
\setcounter{page}{1}

\title[Ky Fan inequalities]{Ky Fan inequalities}

\author[M.S. Moslehian]{Mohammad Sal Moslehian}

\address{Department of Pure Mathematics, Ferdowsi University of
Mashhad, P.O. Box 1159, Mashhad 91775, Iran.}
\email{moslehian@ferdowsi.um.ac.ir, moslehian@member.ams.org}
\urladdr{\url{http://profsite.um.ac.ir/~moslehian/}}

\dedicatory{Dedicated to the memory of Professor Ky Fan}
\subjclass[2010]{Primary 00-02; Secondary 26D15, 47A63.}

\keywords{Ky Fan inequalities, eigenvalue inequality, matrix
inequality, minimax inequality, M-Matrix, Sz\'asz's inequality, Ky
Fan norm, Kantorovich inequality, Ky Fan--Todd determinantal
inequality, Ky Fan--Taussky--Todd inequality.}

\begin{abstract}
There are several inequalities in the literature carrying the name
of Ky Fan. We survey these well-known Ky Fan inequalities and some
other significant inequalities generalized by Ky Fan and review some
of their recent developments.
\end{abstract} \maketitle

\section{Introduction}

According to MathSciNet, Professor Ky Fan (1914-2010) published 126
papers and books and coauthored with 15 mathematicians. His earliest
indexed item goes back to 1940 \cite{1940}. It is notable that
contributions of Ky Fan to mathematics have provided a lot of
influence in the development of nonlinear analysis, convex analysis,
approximation theory, operator theory, linear algebra, mathematical
programming and mathematical economics; see e.g. \cite{LIN}. In the
literature, there are several inequalities due to Ky Fan in various
fields; cf. \cite{BUL}. In this article we try to briefly survey
most important ones by using some information in MathSciNet and
Zentralblatt MATH for some unavailable old papers.

\section{Preliminaries on matrix analysis}

An $n$-tuple $x=(x_1,\ldots,x_n) \in \mathbb{R}^n$ with
$x_1\ge\cdots\geq x_n$ is said to be weakly majorized by an
$n$-tuple $y=(y_1,\ldots,y_n)$ of real numbers with
$y_1\ge\cdots\geq y_n$, denoted by $x \prec_w y$, if
$\sum_{i=1}^kx_i \leq \sum_{i=1}^ky_i$ for $k=1,\ldots,n$. If, in
addition, $\sum_{i=1}^nx_i=\sum_{i=1}^ny_i$, then $x$ is said to be
majorized by $y$ and it is denoted by $x \prec y$. Let
$\mbox{Re}(x)$ stand for $(\mbox{Re}(x_1), \ldots, \mbox{Re}(x_n))$
for an $x=(x_1,\ldots,x_n) \in \mathbb{C}^n$.

Let $\mathbb{M}_n$ be the matrix algebra of all $n \times n$
matrices with entries in the complex field identified with the
algebra  $\mathbb{B}(\mathbb{C}^n)$ of all linear operators on the
Hilbert space $(\mathbb{C}^n, \langle\cdot,\cdot\rangle)$. We denote
by $\mathbb{H}_{n}$ the set of all Hermitian matrices in
$\mathbb{M}_{n}$. By $I_n$ (or $I$ if there is no ambiguity) we
denote the identity matrix of $\mathbb{M}_{n}$. A matrix $A\in
\mathbb{H}_n$ is called positive-semidefinite if $\langle Ax,
x\rangle \geq 0$ holds for every $x\in \mathbb{C}^n$ and then we
write $A\geq 0$. In particular, if $A$ is invertible, then it is
called positive-definite and we write $A>0$. The absolute value of
$A$ is defined by $|A|=(A^*A)^{1/2}$, where $A^*$ denotes the
conjugate transpose of $A$. For $A,B \in \mathbb{H}_n$, we say
$A\leq B$ if $B-A\geq0$.

For a matrix $A\in \mathbb{H}_{n}$, we denote by $\lambda_{1}(A)\geq
\lambda _{2}(A)\geq \cdots \geq \lambda _{n}(A)$ the eigenvalues of
$A$ arranged in the decreasing order with their multiplicities
counted. The notation $\lambda(A)$ stands for the row vector
$(\lambda_{1}(A), \lambda _{2}(A), \ldots, \lambda _{n}(A))$. The
eigenvalue inequality $\lambda(A) \leq \lambda(B)$ means
$\lambda_j(A)\leq \lambda_j(B)\,\,(j=1, 2, \ldots, n)$. The trace
and determinant of a matrix $A \in \mathbb{M}_n$ are defined by
$\mbox{tr}(A)=\sum_{i=1}^n\lambda_i(A)$ and
$\det(A)=\prod_{i=1}^n\lambda_i(A)$, respectively. A norm
$\left\vert \left\vert \left\vert \cdot \right\vert \right\vert
\right\vert $ on $\mathbb{M}_{n}$ is said to be unitarily invariant
if $\left\vert \left\vert \left\vert UAV\right\vert \right\vert
\right\vert =\left\vert \left\vert \left\vert A\right\vert \right\vert \right\vert $ for all $A\in \mathbb{M}_{n}$ and all unitary matrices $%
U,V\in \mathbb{M}_{n}$. The Ky Fan norms, defined as $\| A\|
_{(k)}=\sum_{j=1}^{k}s_{j}(A)$ for $k=1, 2, \ldots, n$, provide a
significant family of unitarily invariant norms; see also \cite{BHA,
MO}. Here $s(A)=(s_{1}(A), \ldots, s_{n}(A))$ denotes the $n$-tuple
of the singular values of $A$, i.e. the eigenvalues of $|A|$,
arranged in the decreasing order.

\section{Ky Fan matrix inequalities}

The contributions of Ky Fan to matrix theory and operator theory are
substantial. In this section we study matrix inequalities due to Ky
Fan. The interested reader is referred to \cite{MO, BHA, HJ} to see
details of proofs.

Ky Fan \cite[1949]{FAN49} proved that if $H, K \in \mathbb{H}_n$,
then $$\lambda(H+K)\prec \lambda(H)+\lambda(K)\qquad (\mbox{Ky Fan
eigenvalue inequality})\,,$$ see also \cite[Theorme 7.14]{ZHA}. An
extension of this result was given by Aujla and Silva \cite{AS} as
$\lambda(f(\alpha A+(1-\alpha)B)) \prec_w \lambda(\alpha
f(A)+(1-\alpha)f(B))$, where $f$ is a convex function on an interval
$J$ and $A, B \in \mathbb{H}_n$ with spectra in $J$. In addition, Ky
Fan \cite[1949]{FAN49} showed that if $A \in \mathbb{M}_n$, $t$ is a
positive integer and the eigenvalues $\kappa_i^{(t)}$ of $(A^t)^\ast
A^t$ are arranged in the decreasing order
$\kappa_1^{(t)}\geq\cdots\geq\kappa_n^{(t)},$ then
$$\kappa_1^{(t)}+\cdots+\kappa_q^{(t)}\leq\kappa_1^t+\cdots+\kappa_q^t\quad
(1\leq q\leq n)\,,$$ where $k_1, \ldots, k_n$ are eigenvalues of
$A^*A$. Ky Fan \cite[1950]{FAN50} established two results which are
powerful tools in obtaining matrix inequalities. They assert that if
$H \in  \mathbb{H}_n$, then
$$\max_{UU^*=I_k}\,\mbox{tr}(UHU^*)=\sum_{i=1}^k\lambda_i(H), \quad \min_{UU^*=I_k}\,\mbox{tr}(UHU^*)=\sum_{i=1}^k\lambda_{n-i+1}(H)\quad (1 \leq k \leq n)$$
where maximum and minimum are taken over all $k \times n$ complex
matrices $U$ satisfying $UU^*=I_k$.

The Ky Fan dominance theorem \cite[1951]{FAN51} states that for $A,
B \in \mathbb{M}_n$, the inequalities $\| A\| _{(k)} \leq \|B\|
_{(k)}\,\,(1 \leq k\leq n)$ hold if and only if $\left\vert
\left\vert \left\vert A\right\vert \right\vert\right\vert \leq
\left\vert \left\vert \left\vert B\right\vert \right\vert
\right\vert $ for all unitarily invariant norms $\left\vert
\left\vert \left\vert \cdot \right\vert \right\vert\right\vert $;
see \cite{MOS2} for an application of the Ky Fan dominance theorem.

Ky Fan \cite[1951]{FAN51} extended an inequality of Weyl \cite{WEY}
(the case $m=1$) and von Neumann \cite{VON} (the case $m=2$): Let
$A_1, \ldots, A_m \in \mathbb{M}_n$ and let $s_{j1}\geq
s_{j2}\geq\cdots \geq s_{jn}$ be the singular values of $A_j$. Then
$$ \max\left|\sum_{i=1}^n\langle U_1A_1\cdots U_mA_mx_i,x_i \rangle\right|\leq\sum_{i=1}^n(s_{1i}\cdots s_{mi}) $$
and
$$ \max\bigg|\det_{1\leq i,k\leq n}\,\left[\langle U_1A_1\cdots U_mA_mx_i,x_k\rangle\right]\bigg| \leq \prod_{i=1}^n (s_{1i}\cdots s_{mi})\,,$$
where $U_1,\ldots, U_m$ run over all unitary matrices and $\{x_1,\ldots, x_n\}$ runs over all orthonormal sets in $\mathbb{C}^n$.\\
Furthermore, Ky Fan showed that for any matrices $A, B \in
\mathbb{M}_n$,
$$s(A+B) \prec_w s(A)+s(B)\qquad (\mbox{Ky Fan singular value inequality}),$$
as well as
$$s_{r+t+1}(A+B) \leq s_{r+1}(A)+s_{t+1}(B),$$
where $t \geq 0, r \geq 0, r+t+1 \leq n$.\\
In the same paper \cite[1951]{FAN51}, it was compared by Ky Fan that
the real parts and modules of the diagonal elements $a_1, \ldots,
a_n$ of an arbitrary matrix $A\in \mathbb{M}_n$ by establishing that
$$\mbox{Re}(a_1, \ldots, a_n) \prec_w (|a_1|, \ldots, |a_n|) \prec_w s(A)\,.$$
Due to the fact that any matrix is unitarily equivalent to an upper
triangular matrix, one gets another Ky Fan inequality
\cite[1950]{FAN50} for eigenvalues as $$\mbox{Re}(\lambda(A)) \prec
\lambda(\mbox{Re}(A)),$$ where $\mbox{Re}(A)=(A+A^*)/2$; cf.
\cite[Lemma 4.20]{ZHX}. A similar result for the eigenvalues of
$\mbox{Im}(A)=(A-A^*)/(2\mbox{i})$ is presented by Amir-Moez and
Horn \cite{AH}. The authors of \cite[1955]{FH} proved that
$$|\lambda_i(\mbox{Re}(A))| \leq s_i(A)\qquad(i=1, \ldots, n)$$ for
any matrix $A \in \mathbb{M}_n$.

Several inequalities for eigenvalues and singular values of two Hermitian matrices $H$, $K$ and $H+{\rm i}K$ is surveyed by Cheng, Horn and Li \cite{CHL}. There are still other inequalities due to Ky Fan.\\
Ky Fan \cite[1954]{FAN54} compared the characteristic roots of
diagonal blocks of a Hermitian matrix and those of the matrix
itself. Ky Fan proved that if $H=\left( \begin{array}{cc} H_{11} &
H_{12}
\\ H_{21}& H_{22} \end{array}\right)$ is an $n \times n$ Hermitian
block matrix and $C(H)=\left( \begin{array}{cc} H_{11} & 0 \\ 0&
H_{22} \end{array}\right)$ is its pinching, then $\lambda(C(H))
\prec \lambda(H)$. This is true for any pinching of $H$; see
\cite[page 50]{BHA}. In the next year, Ky Fan \cite[1955]{FAN55}
proved that for $H$ positive-definite
$$\prod_{i=1}^m\lambda_{n-i+1}(H) \leq
\prod_{i=1}^m\lambda_{k-i+1}(H_{11})\,\,(1 \leq m \leq k),$$ where
the matrix $H_{11}$ and $H_{22}$ are $k \times k$ and $m \times m$
matrices, respectively.

Ky Fan collaborated with Hoffman \cite[1955]{FH} and proved the
following: Let $|||\cdot|||$ be a unitarily invariant norm on
$\mathbb{M}_n$. Then

(i) If $A=U|A|$ is the polar decomposition of $A$, where $U$ is
unitary, and $V$ is any unitary matrix, then $$
|||A-U|||\leq|||A-V|||\leq|||A+U|||\,.$$

(ii) Let $H \in \mathbb{H}_n$. It follows from the identity $A-
\mbox{Re}(A)=\frac{A-H}{2{\rm i}}-\frac{(A-H)^*}{2{\rm i}}$ that $$
|||A-\mbox{Re}(A)|||\leq|||A-H|||\,.$$

(iii) If $H, K \in \mathbb{H}_n$ and $ U=(H-{\rm i}I)(H+{\rm
i}I)^{-1}$, $V=(K-{\rm i}I)(K+{\rm i}I)^{-1}$ are their Cayley
transforms, then
$$ |||U-V|||\leq 2|||H-K|||\,.$$

Ky Fan also worked on systems of inequalities. For instance, Ky Fan
\cite[1963]{BF} gave some necessary and sufficient conditions for
the existence of Hermitian matrices $X_j$ satisfying the system of
linear inequalities $$\sum_{j=1}^n(A_{ij}X_j+X_jA_{ij}{}^\ast)\geq
B_i\,\,(1\leq i\leq m)\quad{\rm and}\quad
\text{tr}\left(\sum_{j=1}^nC_jX_j\right)\geq c\,,$$ where $A_{ij}\in
\mathbb{M}_n$ and  $B_i, C_j\in \mathbb{H}_n$ and $c$ is a real
number.

\section{Results of Ky Fan involving integrals}

For real functions $f$ and $g$ on $[0,1]$, we say that $f \prec g$
whenever $\int_0^xf(t)dt\leq\int_0^xg(t)dt$ for all $0\leq x\leq 1$
and $\int_0^1f(t)dt=\int_0^1g(t)dt$. In \cite[1954]{FL}, Ky Fan and
Lorentz gave sufficient and necessary conditions in terms of mixed
second differences of a suitable differentiable function $\Phi$ in
order that $$
\int_0^1\Phi(t,f_1(t),\ldots,f_n(t))dt\leq\int_0^1\Phi(t,g_1(t),\ldots,g_n(t))dt
$$ for any bounded functions $f_i,g_i$ with $f_i\boldsymbol\prec
g_i$.

\section{Ky Fan--Taussky--Todd inequality}

The Wirtinger inequality reads as follows: If $f$ is a periodic
function of the period $2\pi$ with $f, f' \in L^2([0,2\pi])$ and
$\int_0^{2\pi}f=0$, then
$$\int_0^{2\pi}f^2 \leq \int_0^{2\pi}f'^2$$
with equality if and only if $f(x)=a\cos x+b\sin x$ for some
constants $a, b$. Ky Fan, Taussky and Todd \cite[1955]{FTT} proved
discrete analogs of Wirtinger type. They showed that if $a_1,\dots
,a_n$ are real numbers and $a_0=a_{n+1}=0$, then
$$2\left(1-\cos\frac\pi{n+1}\right)\sum_{k=1}^na_k^2 \leq \sum_{k=1}^{n+1}(a_k-a_{k-1})^2$$
with equality if and only if $a_k=c\sin(k\pi/(n+1))\,\,(k=1, \ldots,
n)$ for some real constant $c$.

They also proved that if $a_1,\dots ,a_n$ are real numbers and
$a_0=0$, then
$$2\left(1-\cos\frac\pi{2n+1}\right)\sum_{k=1}^na_k^2\leq \sum_{k=1}^n(a_k-a_{k-1})^2$$
with equality if and only if $a_k=c\sin(k\pi/(2n+1))\,\,(k=1,
\ldots, n)$, where $c$ is a real constant. Both constants
$2\left(1-\cos\frac\pi{n+1}\right)$ and
$2\left(1-\cos\frac\pi{2n+1}\right)$ are best possible. The
converses of these two inequalities are given by Alzer \cite{ALZ91}.
Another related inequality is that of Ozeki; see \cite{MOS2} for its
operator versions.

\section{Ky Fan--Todd determinantal inequality}

Ky Fan and Todd \cite[1955]{FT} showed that if $\left(
\begin{array}{ccc} a_1 & \ldots &a_n \\ b_1& \ldots & b_n
\end{array}\right)$ is a matrix none of whose two-column minor is
singular and $(p_{ij})$ is a symmetric matrix such that
$p=\sum_{1\leq i<j\leq n}p_{ij} \neq 0$, then
$$\left[\sum_{i=1}^n a_i^2\right]\left[\sum_{i=1}^n a_i^2\sum_{i=1}^n b_i^2-(\sum_{i=1}^n a_ib_i)^2\right]^{-1} \leq p^{-2}\sum_{i=1}^n\left[\sum_{j=1, j\neq i}^np_{ij}a_j(a_jb_i-a_ib_j)^{-1}\right]^2\,. $$
The special case where $p_{ij}= 1$, $a_i=\cos\theta_i$ and
$b_i=\sin\theta_i$ for some $\theta_i$ and all $i, j$ was obtained
by Chassan \cite{CHA} by some statistical arguments. Also a
refinement of it was given by Chong \cite{CHO}.

It is worthy noting that there are a variety of other determinantal
inequalities due to Ky Fan. Ky Fan \cite[1953]{FAN53} (see \cite[p.
214]{MPF}) proved that if $|C|_k$ denotes the product of the first
$k$-th smallest eigenvalues of a real positive-definite matrix $C$
and if $A, B \in \mathbb{M}_n$ are real positive-definite matrices
and $0 \leq \lambda \leq 1$, then

$$|\lambda A + (1-\lambda)B|_k \geq |A|_k^\lambda|B|_k^{1-\lambda}\,.$$

Also Ky Fan \cite[1955]{FAN55PCPS} (see \cite[p. 687]{MO}) proved
that if $H=\left( \begin{array}{cc} H_{11} & H_{12} \\ H_{21}&
H_{22} \end{array}\right)$ and $K=\left( \begin{array}{cc} K_{11} &
K_{12} \\ K_{21}& K_{22} \end{array}\right)$ are positive-definite
$n \times n$ matrices and $H_{11}$ and $K_{11}$ are $k \times k$
matrices, then
$$\left(\frac{\det(H+K)}{\det(H_{11}+K_{11})}\right)^{\frac{1}{n-k}} \geq \left(\frac{\det(H)}{\det(H_{11})}\right)^{\frac{1}{n-k}}+ \left(\frac{\det(K)}{\det(K_{11})}\right)^{\frac{1}{n-k}}\,.$$
In fact, Ky Fan studied many  convex/increasing functions on certain
subsets of Hermitian matrices, e.g.
$f_k(H)=\sum_{i=1}^k\lambda_i(H)\,\,(H \in \mathbb{H}_n)$, from
which one can get some inequalities such as the latter one.

In the same paper \cite[1955]{FAN55PCPS}, Ky Fan proved that if
$(i_1,i_2,\ldots,i_k)$ denotes the principal submatrix of a
positive-definite matrix $H \in \mathbb{H}_n$ with eigenvalues
$\lambda_1\leq\lambda_2\leq\cdots\leq\lambda_n$, then $$
\lambda_1\lambda_2\ldots\lambda_{h+k}\leq
\det(p+1,p+2,\ldots,p+k)\prod_{i=1}^h\frac{\det(i,p+1,p+2,\ldots,n)}{\det(p+1,p+2,\ldots,n)}$$
in which $1\leq h\leq p<n$ and $0\leq k\leq n-p$ and for $k=0$ we
assume that $\det(p+1,\ldots,p+0)=1$.

If $A \in \mathbb{M}_n$ is a real matrix such that $A+A^* >0$, then
Ky Fan \cite[1973]{FAN73} showed that
$$\det(\mbox{Re}(A)) \leq \det(A)\frac{1+\mbox{Re}(\lambda_i(A^{-1}A^*))}{2}$$
for all $i$. If $n \geq 4$ and $A-A^*$ is invertible, then
$$0 < \det\left(\frac{A-A^*}{2}\right) <\det(A)\frac{1-\mbox{Re}(\lambda_i(A^{-1}A^*))}{2}$$
for all $i$.

Ky Fan \cite[1983]{FAN83} used the Jensen inequality to prove that
if $A=H+iK$ with $H$ and $K$ Hermitian is normalizable in the sense
that there is an invertible matrix $T$ such that $T^*AT$ is normal,
and if $p\geq 2/n$, where $n$ is the dimension of the Hilbert space,
then $|\det(A)|^p\geq|\det(H)|^p+|\det(K)|^p$. Ky Fan showed that
equality holds for $p >2/n$ if and only if either $\det(A)=0$ or
$H=0$ (or $k=0$). Ky Fan \cite[1974]{FAN74} already proved the same
result for a strictly dissipative matrix, that is, a matrix whose
imaginary part $\mbox{Im}(A)=(A-A^*)/2\mbox{i}$ is
positive-definite.

\section{Ky Fan mean inequality}

Suppose that $x_1, \ldots, x_n \in (0,\tfrac 12)$,
$\mathcal{A}_n:=(1/n)\sum_{i=1}^nx_i,
\mathcal{G}_n:=(\prod_{i=1}^nx_i)^{1/n},
\mathcal{H}_n:=n/\sum_{i=1}^n(1/x_i)$ and
$\mathcal{A}'_n,\mathcal{G}'_n, \mathcal{H}'_n$ denote the
unweighted arithmetic, geometric and harmonic means of the numbers
$x_i$ and $1-x_i$, respectively. The Ky Fan inequality for means
states that $$\mathcal{G}_n/\mathcal{G}'_n \leq
\mathcal{A}_n/\mathcal{A}'_n \qquad (\mbox{Ky Fan mean
inequality})$$ with equality only if $x_1 =\cdots=x_n$. It first
appeared in the book ``Inequalities'' by Beckenbach and Bellman
\cite{BB}. It is proved by a forward and backward induction as used
for establishing the arithmetic-geometric mean inequality. It can be
also proved by applying Henrici inequality
$\sum^n_{i=1}1/(1+b_i)\geq
n/(1+(\prod^n_{i=1}b_i)^{1/n})\,\,(b_i\geq 1)$ to $b_i= 1/x_i-1\
(i=1,\ldots, n)$. Utilizing the monotonicity of $(1-x^n)/n$ as a
function of $n$, Rooin \cite{ROO} established some Ky Fan type
inequalities. Applying the convexity of the function
$f(x)=(e^{x}-\alpha )/(1+e^{x})\,\,(\alpha>0)$ and the Jensen
inequality, the authors of \cite{DE} generalized the Ky Fan
inequality. Of course, Ky Fan \cite[1951]{FAN51} proved that the
function
$$f(x)=\frac{\left(\prod_{i=1}^n(1-x_i)\right)^{1/n}}{\sum_{i=1}^n(1-x_i)} \cdot \frac{\sum_{i=1}^nx_i}{\left(\prod_{i=1}^nx_i\right)^{1/n}}$$
is symmetric and Schur-convex.

Many extensions, refinements as well as some counterparts have been
given by many authors. Levinson \cite{LEV} extended the inequality
by showing that if $\varphi(u)$ has third derivative for $0<u<2b$,
with $\varphi''(u)\geq 0$, $0<x_i\leq b$ and $0<p_i$,
$i=1,2,\ldots,n$, then
$$\frac{\sum_{i=1}^np_i\varphi(x_i)}{\sum_{i=1}^np_i}-\varphi\left(\frac{\sum_{i=1}^np_ix_i}{\sum_{i=1}^np_i}\right)
\leq\frac{\sum_{i=1}^np_i\varphi(2b-x_i)}{\sum_{i=1}^np_i}-\varphi\left(\frac{\sum_{i=1}^np_i(2b-x_i)}{\sum_{i=1}^np_i}\right).$$
The Ky Fan inequality is the special case when $\varphi(u)=\log u$,
$p_i=1$ and $b=1/2$; see also \cite{DRA}. It is noticed in
\cite{JPS} that Ky Fan's inequality can be deduced from an
inequality proved by L. Lewent in 1908. The interested reader may be
referred to the interesting survey \cite{ALZ2} for some of
generalizations and refinements of the Ky Fan inequality until 1995.

In \cite{ALZ3}, the author presented the following additive version
of the Ky Fan inequality:
$$\min_{1\le i\le n}\frac{x_i}{1-x_i}<\frac{\mathcal{A}_n'-\mathcal{G}_n'}{\mathcal{A}_n-\mathcal{G}_n}<\max_{1\le i\le n}\frac{x_i}{1-x_i}.$$
Alzer, Ruscheweyh and Salinas \cite{ARS} provided the following
refinement of the Ky Fan inequality: $$
\left(\frac{\mathcal{H}_n}{\mathcal{H}'_n}\right)^{n-1}\cdot
\frac{\mathcal{A}_n}{\mathcal{A}'_n}\leq
\left(\frac{\mathcal{G}_n}{\mathcal{G}'_n}\right)^n. $$ Some
refinements of the Ky Fan inequality including means of two or more
variables are obtained by Neuman and Sandor \cite{NS2}. A unified
approach to Sierpi\'nski and Ky Fan inequalities by utilizing a
refined Cauchy inequality and convexity of suitable functions was
presented by Alzer, Ando and Nakamura \cite{AAN}.

\section{Ky Fan inequalities involving $M$-matrices}

A real matrix $A=(a_{ij}) \in \mathbb{M}_n$ is an $M$-matrix if it
can be written as $A=rI-A_1$, where $r$ is greater than the spectral
radius of $A_1$ and $A_1$ has nonnegative entries. This notion was
introduced by Ostrowski \cite{OST}. Ky Fan \cite[1964]{FAN64} gave
several equivalent definitions for the notion of $M$-matrix. One of
them is that $A$ is an $M$-matrix if and only if $A$ is nonsingular
and $A^{-1}$ has nonnegative entries. A complex matrix $B=(b_{ij})$
dominates $A$ if $a_{ii}\leq|b_{ii}|$ for all $i$ and
$|b_{ij}|\leq|a_{ij}|$ for all $i\neq j$. If $B$ is also an
$M$-matrix, then we say that $B$ proportionally dominates $A$
whenever there exist $p_i>0$ such that $p_ia_{ij}\leq b_{ij}$ and
$a_{ij}p_j\leq b_{ij}$ for all $i,j$. If
$\alpha\subseteq\{1,2,\ldots,n\}$, let $A(\alpha)$ denote the
determinant of the principal submatrix formed by the rows and
columns with indices contained in $\alpha$, and we assume that
$A(\emptyset)=1$.

An inequality of Ostrowski \cite{OST} states that
$\det(A)\leq|\det(B)|$ when $B$ dominates $A$, in particular
$\prod_{i=1}^n\lambda_i(A) \leq \prod_{i=1}^nh_{ii}$, when $A$ is an
$M$-matrix. Ky Fan \cite[1966]{FAN66} extended the former inequality
by showing that under conditions above on the entries of $A$ and $B$
if $\alpha_1,\ldots,\alpha_m$ are subsets of $\{1,2,\ldots,n\}$ such
that each of the indices $1,2,\ldots,n$ is contained in at most $k$
of the sets $\alpha_i$, then $$
\det(A)^k/\prod_{i=1}^mA(\alpha_i)\leq|\det(B)^k/\prod_{i=1}^mB(\alpha_i)|.$$

In \cite[1967]{FAN67Shisha} Ky Fan proved that if $A$ and $B$ are
$M$-matrices of order $n$ and a matrix $D\in \mathbb{M}_n$ dominates
$A+B$ and $B$  proportionally dominates $A$, then
$|\det(D)|^{1/n}\geq \det(A)^{1/n} + \det(B)^{1/n}$ and
\begin{eqnarray*}
&&\left|\det(D)/D(\{k+1,\ldots,n\})\right|^{1/k}\\
&&\qquad\qquad\qquad
\geq\left[\det(A)/A(\{k+1,\ldots,n\})\right]^{1/k}+\left[\det(B)/B(\{k+1,\ldots,n\})\right]^{1/k}
\end{eqnarray*}
for all $1\leq k\leq n-1$.

Also, Ky Fan \cite[1960]{FAN60} showed that if $A,B \in\mathbb{M}_n$
are $M$-matrices such that $a_{ij}\leq b_{ij}$ for all $i,j$, and
$$\Phi(A;\alpha, \beta,
\gamma):=\frac{A(\alpha\cap\beta)A(\alpha\cap\gamma)A(\beta\cap\gamma)A(\alpha\cup\beta\cup\gamma)}{A(\alpha)A(\beta)A(\gamma)A(\alpha\cap\beta\cap\gamma)}\,,$$
then $\Phi(A;\alpha, \beta, \gamma) \leq \Phi(A;\alpha, \beta,
\gamma)$ for any subsets $\alpha,\beta,\gamma$ of
$\{1,2,\ldots,n\}$.

Paper \cite[1964]{FAN64} deals with inequalities for principal
subdeterminants of a special product of $M$-matrices. In this paper
Ky Fan introduced the so-called Ky Fan product $A\odot B=
\{\smallmatrix a_{ii}b_{ii},\,\,\,i=j
\\-a_{ij}b_{ij},\,\,\,i\neq j \endsmallmatrix $ and established that for two $M$-matrices $A, B$,
$$1-\Phi(A\odot B;\alpha, \beta, \gamma) \leq (1-\Phi(A;\alpha, \beta, \gamma))(1-\Phi(B;\alpha, \beta, \gamma))\,.$$

For recent studies on some Ky Fan's results about $M$-matrices and
the Ky Fan product see \cite{fprod1}. The importance of the Ky Fan
product is that the set of $M$-matrices is closed under this
product.

\section{Ky Fan generalization of Sz\'asz's inequality}

Suppose that $A \in \mathbb{M}_n$ is a positive-definite matrix and
$P_k \,\,(1\leq k)$ is the product of all its $k\times k$ principal
minors. The Sz\'asz inequality \cite{SZA}, as a generalization of
Hadamard's determinantal inequality, says that
$\prod_{i=1}^na_{ii}=P_1\ge\cdots\ge P_k^{1/( \smallmatrix n-1 \\
k-1
\endsmallmatrix )}\ge\cdots\ge P_n=\det(A)$. Ky Fan
\cite[1967]{FAN67} strengthened the inequality for a certain class
of matrices, named as GKK, including positive-definite matrices,
totally positive matrices and $M$-matrices. By a GKK-matrix we mean
one that all its principal minors are positive and the product of
any two symmetrically situated almost principal minors is real and
nonnegative. His result reads as follows.

Let $A \in \mathbb{M}_n$ be GKK and let $R_0=I$, $R_k=Q_k^{1/(
\smallmatrix p \\ k \endsmallmatrix )},\,(1 \leq k \leq k)$, where
$Q_k=\prod_{i_1<\cdots<i_k}A(\alpha_{i_1}\cup\cdots\cup\alpha_{i_k})$,
in which $\alpha_i\,\,(1\leq k\leq p)$ are $p (\geq 2)$ pairwise
disjoint possibly empty subsets of $\{1,\ldots,n\}$. Then $R_k^2\geq
R_{k-1}R_{k+1}$ and $Q_k^{1/( \smallmatrix p-1 \\ k-1
\endsmallmatrix )}\geq Q_{k+1}^{1/( \smallmatrix p-1 \\ k
\endsmallmatrix )}\,\,(1 \leq k \leq p-1)$, which gives rise to the
Sz\'asz inequality for $p=n$ and $\alpha_i=\{i\}$.

Ky Fan \cite[1992]{FAN92} generalized its former result by removing
the assumption of the pairwise disjointness of the $\alpha_i$, and
showing that it remains true if $A-I$ is either positive-definite or
an $M$-matrix.

\section{Ky Fan extension of Kantorovich inequality}

The Kantorovich inequality asserts that $$\langle Hx,x\rangle\langle
H^{-1}x,x\rangle\leq
(\lambda_1(H)+\lambda_n(H))^2/4\lambda_1(H)\lambda_n(H)\,$$ where
$x$ is a unit vector in $\bold{C}^n$ and $H$ is an $n\times n$
positive-definite matrix. Ky Fan \cite[1966]{FAN66} generalized the
inequality above by showing that if $0 < mI \leq H \leq MI$,
$x_1,\ldots,x_m$ are vectors in $\mathbb{C}^n$ such that
$\sum_{i=1}^n\|x_i\|^2=1$, then
$$\sum_{j=1}^m\langle H^px_j,x_j\rangle\left[\sum_{j=1}^m\langle Hx_j,x_j\rangle\right]^{-p}\leq (p-1)^{p-1}p^{-p}\left(b^p-a^p)^p((b-a)(ab^p-ba^p)^{p-1}\right)^{-1}\,,$$
where $p$ is any integer different from $0$ and $1$. In 1997 Mond
and Pe\v{c}ari\'c \cite{MP} gave an operator version of Ky Fan's
inequality above. Another extension of Kantorovich inequality was
given by Furuta \cite{FUR}. Ky Fan proved that if $H, K$ are
positive operators on a Hilbert space, $H\ge K>0$ and $MI\ge K\ge
mI>0$, then $$\left(\frac{M}{m}\right)^{p-1}H^p\ge
\frac{(p-1)^{p-1}}{p^p}\frac{(M^p-m^p)^p}{(M-m)(mM^p-Mm^p)^{p-1}}H^p\ge
K^p$$ holds for all $p\geq 1$. The constant
$\kappa_+(m,M,p)=\frac{(p-1)^{p-1}}{p^p}\frac{(M^p-m^p)^p}{(M-m)(mM^p-Mm^p)^{p-1}}$
is called the Ky Fan--Furuta constant in the literature; cf.
\cite{fmps}.

\section{Ky Fan minimax inequality}

In his seminal work \cite[1972]{FAN72}, Ky Fan presented the
following result:

Let $S$ be a nonempty compact convex subset of a Hausdorff
topological vector space $X$, and let $f(x,y)$ be a function from $S
\times S$ to the real numbers that is a lower semicontinuous
function of $y$ for each fixed $x \in S$ and is a quasi-concave
function of $x$ for each fixed $y \in S$. Then the following minimax
inequality holds:
$$\displaystyle{\min_{y\in S}\max_{x\in S} f(x,y) \leq \sup_{x \in S}f(x,x)} \qquad (\mbox{Ky Fan minimax inequality}).$$

It is equivalent to the Brouwer--Kakutani fixed point theorem and
plays an essential role in several fields such as game theory,
mathematical economics, variational inequalities, fixed point
theory, control theory as well as nonlinear and convex analysis.

This inequality has evoked the interest of many mathematicians:
Ferro \cite{FER} got a minimax inequality for vector-valued
functions. Ding \cite{DIN} obtained a minimax inequality by a
generalized relatively Knaster--Kuratowski--Mazurkiewicz (R-KKM)
mapping. Kuwano, Tanaka and Yamada \cite{KTY} obtained a minimax
inequality by using a unified scalarization. Also Hussain, Khan and
Agarwal \cite{HKA} used common fixed point theorem for a condensing
map to obtain some Ky Fan type approximation theorems.

\section{Sharpening of the von Neumann inequality}

The well-known von Neumann inequality asserts that if $A$ is a
proper contraction (i.e, $\|A\|< 1$) acting on a Hilbert space, then
$\|f(A)\|\leq 1$, where $f(A)$ is defined by applying the
Riesz--Dunford functional calculus to a suitable analytic function
$f$. In several papers Ky Fan sharpened von Neumann's inequality by
finding a bound of modulus less than $1$ for $\|f(A)\|$; see
\cite[1987]{FAN872}.

(i) Ky Fan \cite[1979]{FAN791} proved the following inequality
related to Julia's Lemma (for other results see \cite[1979]{AF}):
\begin{align*}
\left[I-\overline{f(z)}f(A)\right]\left[I-f(A)^*f(A)\right]&^{-1}\left[I-f(z)f(A)^*\right]\\&
\leq \frac{1-|f(z)|^2}{1-|z|^2}(I-\bar{z}A)(I-A^*A)^{-1}(I-zA^*)\,.
\end{align*}

(ii) Ky Fan \cite[1987]{FAN871} established that
$$\|f(A)\|\leq[\|\mu_W(A)\|+|f(w)|]\,[1+ \|\mu_W(A)\|\,|f(w)|]^{-1},$$ in which $|w|<1$ and $\mu_W$ is the M\"obius transformation $\mu_W(z)=\frac{z-w}{1-\overline wz}$. \\
For recent develeopments concerning a multivariable von Neumann
inequality see \cite{von}.

\section{Operator Harnack inequalities}

Ky Fan \cite[1988]{FAN88} provided two operator Harnack inequalities:\\
If $F(z)$ is a Hilbert space operator-valued analytic function on
the open unit disk such that  $F(0)=I$  and $\mbox{Re}(F(z))>0$,
then $$(1-|z|)/(1+|z|) I \leq\mbox{Re}(F(z))\leq(1+|z|)/ (1-|z|) I$$
and
$$-2|z|/(1-|z|^2) I \leq\mbox{Im}(F(z))\leq 2|z|/(1- |z|^2) I.$$
See also similar results of Ky Fan in \cite[1980]{FAN80}.\\

\textbf{Acknowledgement.} The author would like to thank Dr. Farzad
Dadipour and anonymous referees for carefully reading the manuscript
and useful comments.

\bibliographystyle{amsplain}

\begin{thebibliography}{99}

\bibitem{ALZ91} H. Alzer, \textit{Converses of two inequalities of Ky Fan, O. Taussky and J. Todd}, J. Math. Anal. Appl. \textbf{161} (1991), 142–-147.

\bibitem{ALZ2} H. Alzer, \textit{The inequality of Ky Fan and related results}, Acta Appl. Math. \textbf{38} (1995), 305--354.

\bibitem{ALZ3} H. Alzer, \textit{On an additive analogue of Ky Fan's inequality}, Indag. Math. (N.S.) \textbf{8} (1997), no. 1, 1--6.

\bibitem{AAN} H. Alzer, T. Ando and Y. Nakamura, \textit{The inequalities of W. Sierpi\'nski and Ky Fan}, J. Math. Anal. Appl. \textbf{149} (1990), no. 2, 497--512.

\bibitem{ARS} H. Alzer, S. Ruscheweyh and L. Salinas, \textit{On Ky Fan-type inequalities}, Aequationes Math. \textbf{62} (2001), no. 3, 310--320.

\bibitem{AH} A.R. Amir-Mo\'ez and A. Horn, \textit{Singular values of a matrix}, Amer. Math. Monthly \textbf{65} (1958), 742–-748.

\bibitem{AF} T. Ando and K. Fan, \textit{Pick--Julia theorems for operators}, Math. Z. \textbf{168} (1979), no. 1, 23–-34.

\bibitem{AS} J.S. Aujla and F.C. Silva, \textit{Weak majorization inequalities and convex functions}, Linear Algebra Appl. \textbf{369} (2003), 217–-233.

\bibitem{BB} E.F. Beckenbach and R. Bellman, \textit{Inequalities}, Ergebnisse der Mathematik und ihrer Grenzgebiete, N. F., Bd. 30 Springer-Verlag, Berlin-G\"ottingen-Heidelberg 1961.

\bibitem{BF} R. Bellman and
 K. Fan, \textit{On systems of linear inequalities in Hermitian matrix variables}, in: V.L. Klee (Ed.), Convexity, The Proceedings of Symposia in Pure Mathematics, Vol. 7, American Mathematical Society, Providence, RI, 1963, pp. 1-–11.

\bibitem{BHA} R. Bhatia, \textit{Matrix Analysis}, Graduate Texts in Mathematics, 169. Springer-Verlag, New York, 1997.

\bibitem{BUL} P.S. Bullen, \textit{A dictionary of inequalities}, Pitman Monographs and Surveys in Pure and Applied Mathematics, 97. Longman, Harlow, 1998.

\bibitem{CHA} J.B. Chassan, \textit{A statistical derivation of a pair of trigonometric inequalities},  Amer. Math. Monthly  \textbf{62} (1955), 353–-356.

\bibitem{CHL} C.-M. Cheng, R.A. Horn and C.-K. Li, \textit{Inequalities and equalities for the Cartesian decomposition of complex matrices}, Linear Algebra Appl. \textbf{341} (2002), 219-–237.

\bibitem{CHO} K.K. Chong, \textit{Refinements of the Fan-Todd's inequalities}, Chinese Ann. Math. Ser. B \textbf{23} (2002), no. 1, 75–-84.

\bibitem{DIN} X.P. Ding, \textit{New generalized R-KKM type theorems in general topological spaces and applications}, Acta Math. Sin. (Engl. Ser.) \textbf{23} (2007), no. 10, 1869--1880.

\bibitem{DE} J.L. Diaz-Barrero, J.J. Egozcue and J. Gibergans-B\'aguena, \textit{On the Ky Fan inequality and some of its applications}, Comput. Math. Appl. \textbf{56} (2008), no. 9, 2279--2284.

\bibitem{DRA} S.S. Dragomir, \textit{Some refinements of Ky Fan's inequality}, J. Math. Anal. Appl. \textbf{163} (1992), no. 2, 317--321.

\bibitem{1940} K. Fan, \textit{Sur les types homogènes de dimensions} (French), C. R. Acad. Sci. Paris  \textbf{211} (1940), 175–-177.

\bibitem{FAN49} K. Fan, \textit{On a theorem of Weyl concerning eigenvalues of linear transformatioins I.},
Proc. Nat. Acad. Sci. U.S.A. \textbf{35} (1949), 652--655.

\bibitem{FAN50} K. Fan, \textit{On a theorem of Weyl concerning eigenvalues of linear transformations. II}, Proc. Nat. Acad. Sci. U. S. A. \textbf{36} (1950), 31-–35.

\bibitem{FAN51} K. Fan, \textit{Maximum properties and inequalities for the eigenvalues of completely continuous operators}, Proc. Nat. Acad. Sci. U.S.A. \textbf{37} (1951), 760--766.

\bibitem{FAN53} K. Fan, \textit{Advanced Problems and Solutions: Solutions: 4430}, Amer. Math. Monthly  \textbf{60} (1953),  no. 1, 50.

\bibitem{FAN54} K. Fan, \textit{Inequalities for eigenvalues of Hermitian matrices}, Contributions to the solution of systems of linear equations and the determination of eigenvalues, pp. 131–-139. National Bureau of Standards Applied Mathematics Series No. 39. U. S. Government Printing Office, Washington, D.C., 1954.

\bibitem{FAN55} K. Fan, \textit{A comparison theorem for eigenvalues of normal matrices}, Pacific J. Math. \textbf{5} (1955), 911–-913.

\bibitem{FAN55PCPS} K. Fan, \textit{Some inequalities concerning positive-definite Hermitian matrices}, Proc. Cambridge Philos. Soc. \textbf{51} (1955), 414–-421.

\bibitem{FAN60} K. Fan, \textit{Note on $M$-matrices}, Quart. J. Math. Oxford Ser. (2) \textbf{11} 1960, 43–-49.

\bibitem{FAN64} K. Fan, \textit{Inequalities for $M$-matrices}, Nederl. Akad. Wetensch. Proc. Ser. A \textbf{67} (= Indag. Math. \textbf{26}) (1964), 602–-610.

\bibitem{FAN66} K. Fan, \textit{Some matrix inequalities}, Abh. Math. Sem. Univ., Hamburg \textbf{29} (1966), 185–-196.

\bibitem{FAN67Shisha} K. Fan, \textit{Inequalities for the sum of two matrices}, in Inequalities, O. Shisha, ed., Academic Press, New York, 1967.

\bibitem{FAN67} K. Fan, \textit{Subadditive functions on a distributive lattice and an extension of Sz\'asz's inequality},
J. Math. Anal. Appl. \textbf{18} (1967), 262-–268.

\bibitem{FAN72} K. Fan, \textit{A minimax inequality and applications}, Inequalities, III (Proc. Third Sympos., Univ. California, Los Angeles, Calif., 1969; dedicated to the memory of Theodore S. Motzkin), pp. 103--113. Academic Press, New York, 1972.

\bibitem{FAN73} K. Fan, \textit{On real matrices with positive-definite symmetric component}, Linear and Multilinear Algebra \textbf{1} (1973), no. 1, 1-–4.

\bibitem{FAN74} K. Fan, \textit{On strictly dissipative matrices}, Linear Algebra and Appl. \textbf{9} (1974), 223–-241.

\bibitem{FAN791} K. Fan, \textit{Julia's lemma for operators}, Math. Ann. \textbf{239} (1979), no. 3, 241-–245.

\bibitem{FAN80} K. Fan, \textit{Harnack's inequalities for operators}, General inequalities, 2 (Proc. Second Internat. Conf., Oberwolfach, 1978), pp. 333–-339, Birkh\"auser, Basel-Boston, Mass., 1980.

\bibitem{FAN83} K. Fan, \textit{Normalizable operators}, Linear Algebra Appl. \textbf{52/53} (1983), 253-–263.

\bibitem{FAN871} K. Fan, \textit{Sharpened forms of an inequality of von Neumann}, Math. Z. \textbf{194} (1987), no. 1, 7–-13.

\bibitem{FAN872} K. Fan, \textit{Applications and sharpened forms of an inequality of von Neumann}, Current trends in matrix theory (Auburn, Ala., 1986), 113–-121, North-Holland, New York, 1987.

\bibitem{FAN88} K. Fan, \textit{Inequalities for proper contractions and strictly dissipative operators}, Linear Algebra Appl. \textbf{105} (1988), 237-–248.

\bibitem{FAN92} K. Fan, \textit{Some inequalities for matrices $A$ such that $A-I$ is positive-definite or an $M$-matrix}, Linear and Multilinear Algebra \textbf{32} (1992), no. 2, 89–-92.

\bibitem{FH} K. Fan and A.J. Hoffman, \textit{Some metric inequalities in the space of matrices}, Proc. Amer.Math. Soc. \textbf{6} (1955) 111–-116.

\bibitem{FL} K. Fan and G.G. Lorentz, \textit{An integral inequality}, Amer. Math. Monthly \textbf{61} (1954), 626-–631.

\bibitem{FTT} K. Fan, O. Taussky and J. Todd, \textit{Discrete analogs of inequalities of Wirtinger}, Monatsh. Math. \textbf{59} (1955), 73-–90.

\bibitem{FT} K. Fan and J. Todd, \textit{A determinantal inequality}, J. London. Math. Soc. \textbf{30} (1955), 58–-64.

\bibitem{FER} F. Ferro, \textit{A minimax theorem for vector-valued functions. II}, J. Optim. Theory Appl. \textbf{68} (1991), no. 1, 35–-48.

\bibitem{FUR} T. Furuta, \textit{Operator inequalities associated with H\"older-McCarthy and Kantorovich inequalities}, J. Inequal. Appl. \textbf{2}  (1998), no. 2, 137--148.

\bibitem{fmps} T. Furuta, J. Mi\'ci\'c Hot, J.E. Pe\v cari\'c and Y. Seo, \textit{Mond--Pecaric Method in Ooperator Inequalities, Inequalities for
Bounded Selfadjoint Operators on a Hilbert Space}, Monographs in
Inequalities 1. Zagreb: Element, 2005.

\bibitem{von} A. Grinshpan, D. Kaliuzhnyi-Verbovetskyi, V. Vinnikov, and H.J. Woerdeman, \textit{Classes of tuples of commuting contractions satisfying the multivariable von Neumann inequality}, J. Funct. Anal. \textbf{256}, 3035-–3054.

\bibitem{HJ} R.A. Horn and C.R Johnson, \textit{Topics in matrix analysis}, Corrected reprint of the 1991 original. Cambridge University Press, Cambridge, 1994.

\bibitem{HKA} N. Hussain A.R. Khan and R.P. Agarwal, \textit{Krasnosel\'skii and Ky Fan type fixed point theorems in ordered Banach spaces},  J. Nonlinear Convex Anal. \textbf{11}  (2010),  no. 3, 475–-489.

\bibitem{JPS} M.V. Jovanovi\'c,  T.K. Pog\'any and J. S\'andor, \textit{Notes on certain inequalities by H\"older, Lewent and Ky Fan}, J. Math. Inequal. \textbf{1} (2007), no. 1, 53--55.

\bibitem{KTY} I. Kuwano, T. Tanaka and S. Yamada, \textit{Unified scalarization for sets and set-valued Ky Fan minimax inequality}, J. Nonlinear Convex Anal. \textbf{11} (2010), no. 3, 513--525.

\bibitem{LEV} N. Levinson, \textit{Generalization of an inequality of Ky Fan}, J. Math. Anal. Appl. \textbf{8} (1964), 133--134.

\bibitem{LIN} B.-L. Lin, \textit{Every waking moment Ky Fan (1914--2010)},  Notices Amer. Math. Soc. \textbf{57} (2010),  no. 11, 1444–-1447.

\bibitem{MO} A.W. Marshall, I. Olkin and B.C. Arnold, \textit{Inequalities: theory of majorization and its applications}, Second edition. Springer Series in Statistics. Springer, New York, 2011.

\bibitem{MOS2} J.S. Matharu, M.S. Moslehian and J.S. Aujla, \textit{Eigenvalue extensions of Bohr's inequality}, Linear Algebra Appl. \textbf{435} (2011), no. 2,  270--276.

\bibitem{MPF} D.S. Mitrinovi\'c, J.E. Pe\v{c}ari\'c and A.M. Fink,\textit{Classical and New Inequalities in Analysis}, Mathematics and its Applications (East European Series), 61. Kluwer Academic Publishers Group, Dordrecht, 1993.

\bibitem{MP} B. Mond and J.E. Pe\v{c}ari\'c, \textit{Matrix inequalities for convex functions}, J. Math. Anal. Appl. \textbf{209} (1997), no. 1, 147–-153.

\bibitem{MOS1} M.S. Moslehian, R. Nakamoto and Y. Seo, \textit{A Diaz-Metcalf type inequality for positive linear maps and its applications}, Electron. J. Linear Algebra \textbf{22} (2011), 179--190.

\bibitem{NS2} E. Neuman and J. Sandor, \textit{On the Ky Fan inequality and related inequalities. II}, Bull. Austral. Math. Soc. \textbf{72} (2005), no. 1, 87--107.

\bibitem{OST} A. Ostrowski, \textit{\"Uber die determinanten mit \"uberwiegender Hauptdiagonale} (German),  Comment. Math. Helv.  \textbf{10}  (1937),  no. 1, 69–-96.

\bibitem{ROO} J. Rooin, \textit{An approach to Ky Fan type inequalities from binomial expansions}, Math. Inequal. Appl. \textbf{11} (2008), no. 4, 679--688.

\bibitem{SZA} O. Sz\'asz, \textit{Über eine Verallgemeinerung des Hadamardschen Determinantensatzes} (German), Monatsh. Math. Phys. \textbf{28} (1917), no. 1, 253–-257.

\bibitem{VON} J. von Neumann,  \textit{Some matrix-inequalities and metrization of matrix-space}, Tomsk. Univ. Rev. 1, (1937), 286–-300. [John von Neumann Collected Works (A. H. Taub, ed.), Vol. IV, pp. 205–218. Pergamon, Oxford, 1962].

\bibitem{fprod1} B.-Y. Wang, X. Zhang and F. Zhang, \textit{On the Hadamard product of inverse $M$-matrices}, Linear Algebra Appl. \textbf{305} (2000), no. 1-3, 23–-31.

\bibitem{WEY} H. Weyl, \textit{Inequalities between the two kinds of eigenvalues of a linear transformation}, Proc. Nat. Acad. Sci. U. S. A. \textbf{35} (1949), 408–-411.

\bibitem{ZHX} X. Zhan, \textit{Matrix Inequalities}, Lecture Notes in Mathematics, 1790. Springer-Verlag, Berlin, 2002.

\bibitem{ZHA} F. Zhang, \textit{Matrix Theory, Basic Results and Techniques}, Universitext, Springer-Verlag, New York, 1999.


\end{thebibliography}

\end{document}